\documentclass[12pt]{article}
\usepackage{amsfonts,euscript}

\renewcommand{\thefootnote}{*,\S}

\newtheorem{thm}{Theorem}[section]

\newtheorem{rem}[thm]{Remark}
\newtheorem{defn}[thm]{Definition}
\newtheorem{conj}[thm]{Conjecture}
\newtheorem{prop}[thm]{Proposition}
\newtheorem{quest}[thm]{Question}
\newtheorem{lem}[thm]{Lemma}
\newtheorem{cor}[thm]{Corollary}
\def\eqref#1{(\ref{#1})}
\def\qed{~\vrule height8pt width 5pt depth -1pt\medskip}
\long\def\symbolfootnote[#1]#2{\begingroup%
\def\thefootnote{\fnsymbol{footnote}}\footnote[#1]{#2}\endgroup}

\begin{document}
\bibliographystyle{plain}

\begin{center}
{\bf \large The prime spectrum of algebras of quadratic growth}\\
\vspace{7 mm} \centerline{Jason P. Bell\symbolfootnote[1]{The first author thanks NSERC for its generous support.}}
\centerline{Department of Mathematics} \centerline{Simon Fraser University}
\centerline{
8888 University Drive} 
\centerline{Burnaby, BC, Canada}
\centerline{V5A 1S6}
\centerline{\tt jpb@math.sfu.ca}
\vskip 7mm
\centerline{ Agata Smoktunowicz\symbolfootnote[2]{The second author was supported by Grant No. EPSRC
 EP/D071674/1.}}
\centerline{Maxwell Institute of Sciences}
\centerline{School of Mathematics, University of Edinburgh}
\centerline{James Clerk Maxwell Building, King's Buildings}
\centerline{Mayfield Road, Edinburgh EH9 3JZ, Scotland}
\centerline{\tt A.Smoktunowicz@ed.ac.uk}

\center{Keywords: GK dimension, quadratic growth, primitive rings, PI rings, graded algebra.}
\end{center}

\begin{abstract}  We study prime algebras of quadratic growth.  Our first result is that if $A$ is a prime monomial algebra of quadratic growth then $A$ has finitely many prime ideals $P$ such that $A/P$ has GK dimension one.  This shows that prime monomial algebras of quadratic growth have bounded matrix images.  We next show that a prime graded algebra of quadratic growth has the property that the intersection of the nonzero prime ideals $P$ such that $A/P$ has GK dimension $2$ is non-empty, provided there is at least one such ideal.  From this we conclude that a prime monomial algebra of quadratic growth is either primitive or has nonzero locally nilpotent Jacobson radical.  Finally, we show that there exists a prime monomial algebra $A$ of GK dimension two with unbounded matrix images and thus the quadratic growth hypothesis is necessary to conclude that there are only finitely many prime ideals such that $A/P$ has GK dimension $1$.

\end{abstract}

\section{Introduction}
We study finitely generated algebras of quadratic growth.  Given a field $k$ and a finitely generated $k$-algebra $A$, a $k$-subspace $V$ of $A$ is called a \emph{frame} of $A$ if $V$ is finite dimensional, $1\in V$, and $V$ generates $A$ as a $k$-algebra.
We say that $A$ has \emph{quadratic growth} if there exist a frame $V$ of $A$ and constants $C_1,C_2>0$ such that
$$C_1 n^2 \ \le \ {\rm dim}_k(V^n) \ \le \ C_2n^2\qquad {\rm for~all}~n\ge 1.$$
We note that an algebra of quadratic growth has GK dimension $2$.  More generally, the GK dimension of a finitely generated $k$-algebra $A$ is defined to be
$${\rm GKdim}(A) \ = \ \limsup_{n\rightarrow\infty} \log\bigg({\rm dim}(V^n)\bigg)\bigg/\log\, n,$$
where $V$ is a frame of $A$.  While algebras of quadratic growth have GK dimension $2$, it is not the case that an algebra of GK dimension $2$ necessarily has quadratic growth.  We note that in 
the case that A is a finitely generated commutative algebra, GK dimension 
is equal to Krull dimension. For this reason, GK dimension has seen great 
use over the years as a useful tool for obtaining noncommutative analogues 
of results from classical algebraic geometry. For more information about GK 
dimension we refer the reader to Krause and Lenagan \cite{KL}.

We consider prime monomial algebras of quadratic growth.  A $k$-algebra $A$ is a \emph{monomial algebra} if 
$$A\ \cong \ k\{x_1,\ldots ,x_d\}/I,$$ for some ideal $I$ generated by monomials in $x_1,\ldots ,x_d$.  Monomial algebras are useful for many reasons.  First, Gr\"obner bases a monomial algebra to a finitely generated algebra, and thus monomial algebras can be used to answer questions about ideal membership and Hilbert series for general algebas.  Second, many questions for algebras reduce to combinatorial problems for monomial algebras and can be studied in terms of forbidden subwords.  For these reasons, monomial algebras are a rich area of study.  The paper of Belov,  Borisenko, and Latyshev \cite{Belov}  gives many useful and interesting results about monomial algebras.

Our first result is the following theorem.  
\begin{thm} \label{thm: main} Let $k$ be a field and let $A$ be a prime monomial $k$-algebra of quadratic growth.  Then the set of primes $P$ such that
${\rm GKdim}(A/P)=1$ is finite; moreover, all such primes are monomial ideals.  In particular, $A$ has bounded matrix images.
\end{thm}
We also show that if there is a constant $C>0$ such that the number of words of length at most $n$ with nonzero image in $A$ is at most $Cn^2$ for $n$ sufficiently large, then $A$ has at most $2C$ primes such that $A/P$ has GK dimension one.  In general, if a finitely generated $k$-algebra $A$ has quadratic growth, we define the \emph{growth constant} of $A$ to be
$${\rm GC}(A) \ := \ \inf_V\, \limsup_{n\rightarrow\infty} {\rm dim}(V^n)/n^2,$$
where the infimum is taken over all frames $V$ of $A$. 
Theorem \ref{thm: main} leads us to make the following conjecture.
\begin{conj} There exists a function $F:(0,\infty)\rightarrow (0,\infty)$ such that for every field $k$ and every finitely generated non-PI prime Noetherian $k$-algebra $A$ of quadratic growth,
$A$ has at most $F({\rm GC}(A))$ prime ideals $P$ such that ${\rm GKdim}(A/P)=1$.
\end{conj}
We know of no example of a finitely generated prime Noetherian $k$-algebra of quadratic growth with more than $2{\rm GC}(A)+1$ primes of co-GK dimension $1$.  

We next turn our attention to graded algebras of quadratic growth.  In this case, we give an analogue of Bergman's gap theorem---which states that there are no algebras of GK dimension strictly between $1$ and $2$---for ideals.
 \begin{thm}\label{thm: asmain1}
 Let $K$ be a field and let $A=\bigoplus_{i=0}^{\infty }A_{i}$ be a
 finitely generated graded $K$-algebra with quadratic growth, 
  generated in degree one.
 Suppose $u$ is a nonzero homogeneous element of $A$ and let $(u)$ denote the two-sided ideal of $A$ generated by $u$. 
Then either:
\begin{enumerate}
\item there is a natural number $m$ such
 that \[{\rm dim}_{K}\left((u) \bigcap \bigoplus_{i=0}^{n}A_{i}\right ) \ \ge \  {(n-m)(n-m-1)\over 2}\]
  for all $n$ sufficiently large; or
  \item  there is  a positive constant $C$ such that
  \[{\rm dim}_{K}\left((u)\bigcap \bigoplus_{i=1}^n A_{i}\right )\ < \ Cn\] for all $n$ sufficiently large.
  \end{enumerate}
   Moreover, if $A$ is prime then the former holds.
 \end{thm}
Using this theorem we are able to obtain the following result about prime ideals in graded algebras of quadratic growth.
\begin{thm} Let $K$ be a field, and let $A=\bigoplus_{i=0}^{\infty }A_{i}$ be a prime affine graded non-PI
 $K$-algebra, generated in degree one.
If $A$ has quadratic growth then the intersection of all nonzero prime ideals $P$ such that $A/P$ has GK dimension $2$ is nonzero, where we take an empty intersection to be all of $A$. \label{thm: asmain2}
\end{thm}
As a result of this theorem and Theorem \ref{thm: main} we obtain the following corollary.

\begin{cor} \label{thm: asmain3} Let $A$ be a finitely generated prime monomial algebra of quadratic growth.  Then $A$ has bounded matrix images and either $A$ is primitive or has nonzero locally nilpotent Jacobson radical.
\end{cor}

We note that a finitely generated prime monomial algebra that has GK dimension greater than $1$ cannot be PI \cite{Belov}.  A conjecture of Small \cite[Question 3.2]{Bell} is that a finitely generated prime Noetherian algebra of quadratic growth is either primitive or PI.  By the Jacobson density theorem, primitive algebras are dense subrings of endomorphism rings over a division algebra.   For this reason, primitive ideals are an important object of study and classifying primitive ideals is often an important intermediate step in classifying finite dimensional representations of an algebra.  
 
Finally, we turn to the more general setting of algebras of GK dimension $2$.  In this case, we are able to show that prime algebras can have unbounded matrix images.

\begin{thm} \label{thm: main2}
There exists a finitely generated prime monomial algebra $A$ of GK dimension $2$ with unbounded matrix images.
\end{thm}
This answers a question of Small \cite[Question 3.1]{Bell}; Small's question is still open when the additional hypothesis that the algebra be Noetherian is added.
Irving \cite{Irving} obtained a similar example of an algebra of GK dimension $2$ with unbounded matrix images; his example is not prime, however.

\section{Prime ideals of co-GK $1$}
We begin our investigation of monomial algebras by giving a description of monomial algebras of GK dimension $1$.  We first introduce the following notation.  Given a finite alphabet $\{x_1,\ldots ,x_d\}$, we call an infinite product of the form $$x_{i_1}x_{i_2}x_{i_3}\cdots $$
a \emph{right infinite word} on $\{x_1,\ldots ,x_d\}$.  Left infinite words are defined analogously and doubly infinite words are both left and right infinite.  
Given a finite alphabet $\{x_1,\ldots ,x_d\}$
 and a right infinite word $W$, we let 
$A_W$ denote the algebra $k\{x_1,\ldots ,x_d\}/I$, where $I$ is the monomial ideal generated by the collection of words which do not occur as a subword of $W$.  In particular, the images of the distinct subwords of $W$ form a basis for $A_W$.
\begin{rem} Let $A$ be a finitely generate prime monomial algebra.  Then there is a right infinite word $W$ on a finite alphabet such that $A\cong A_W$. \label{rem: 1}
\end{rem}
The reason for this is that if $A=k\{x_1,\ldots ,x_d\}/I$, then we can pick an enumeration $W_1,W_2,\ldots $ of the words on $\{x_1,\ldots , x_d\}$ with nonzero image in $A$.  Then since $A$ is prime, there exist words $V_1,V_2,\ldots $ such that $W_1V_1W_2V_2\cdots V_mW_{m+1}$ has nonzero image in $A$ for every $m$.  Let $$W=W_1V_1W_2V_2\cdots .$$
Then $W$ is a right infinite word.  Moreover every subword of $W$ has nonzero image in $A$ since every subword of $W$ is a subword of $W_1V_1W_2V_2\cdots V_mW_{m+1}$ for some $m$.  Moreover, every word with nonzero image in $A$ is a subword of the form $W_i$ for some $i$ and hence is a subword of $W$.  It follows that $A\cong A_W$.

\begin{lem} Let $k$ be a field and let $A$ be a finitely generated prime monomial $k$-algebra of GK dimension $1$.  Then there is a word $W$ on a finite alphabet such that $A\cong A_{W^{\omega}}$, where $W^{\omega}$ is the right infinite word $WWW\cdots .$  Conversely, if $W$ is a word on a finite alphabet, then the algebra $A_{W^{\omega}}$ is a finitely generated prime algebra of GK dimension $1$. \label{lem: 0}
\end{lem}
{\bf Proof.} Since $A$ is prime, there is a right infinite word $U$ such that $A\cong A_U$ by Remark \ref{rem: 1}.
Since $A$ has GK dimension $1$ and the images of the subwords of $U$ form a basis for $A_U$, we see that $U=W_1W^{\omega}$ for some words $W_1$ and $W$ by a result of Bergman \cite[Lemma 2.4]{KL}.  We claim that $W_1W^{\omega}$ is itself a subword of $W^{\omega}$.  To see this, suppose that there is some subword $V$ of $W_1W^{\omega}$ that does not appear as a subword of $W^{\omega}$.  Let $d$ denote the length of $W_1$.  Since $A_U$ is prime we see there exist words $V_1, V_2,\ldots $ such that $VV_1VV_2V\cdots V_dV$ is a subword of $W_1W^{\omega}$.  Moreover since $W_1$ has length $d$, we see that the last occurrence of $V$ in $VV_1VV_2V\cdots V_dV$ must be contained entirely in $W^{\omega}$.  This is a contradiction, we conclude that $A_U=A_{W^{\omega}}$.  The converse follows easily from the fact that $W^{\omega}$ has ${\rm O}(1)$ subwords of length $n$ and each of these words occur infinitely often.  The result follows. \qed

\begin{lem} Let $k$ be a field and let $A=k\{x_1,\ldots ,x_d\}/I$ be a prime monomial algebra of GK dimension at least
$2$.  Let $W$ be a word on the alphabet $\{x_1,\ldots ,x_d\}$ such that every power of $W$ has nonzero image in $A$ and let $n$ be a positive
integer.  Then there exists a right infinite word $U$ on $\{x_1,\ldots ,x_d\}$ such that
every subword of $W^nU$ has nonzero image in $A$ and $U$ does not contain 
$W$ as an initial subword. \label{lem: 1}
\end{lem}
{\bf Proof.} Pick a word $V_1$ that is not a subword of
$W^{\omega}$.  Since $A$ is prime there exists a word
$V_2$ such that 
$W^nV_2V_1$ has nonzero image in $A$.  Since $A$ is prime there
exists an infinite word $V$ with initial subword $V_2V_1$ such that
each subword of $W^nV$ has nonzero image in $A$.  Since $V_1$ is not
subword of $W^{\omega}$, there is some $m\ge n$ such that
$$W^nV=W^mU$$ where $U$ is a right infinite word that does not
have $W$ as an initial subword.  This completes the proof.
\qed  
\vskip 2mm
We first show that prime monomial algebras $A$ of quadratic growth only have finitely many prime monomial ideals $P$ such that $A/P$ has GK dimension $1$. 
\begin{prop} Let $A$ be a prime monomial algebra of quadratic growth. \label{prop: 1}
Then
$A$ has only finitely many prime monomial ideals $P$ such that $A/P$ has GK dimension $1$.
\end{prop}
{\bf Proof.}  We may assume that $A=k\{x_1,\ldots ,x_d\}/I$.  By assumption there is a constant $C>0$ such that there are at most $Cn^2$ words of length at most $n$ with nonzero image in $A$ for every $n\ge 1$.
Let $$\mathcal{S} = \{P~|~{\rm GKdim}(A/P)=1,~ P~{\rm a~monomial~ideal}\}.$$  We show that
$\mathcal{S}$ is finite.  By Lemma \ref{lem: 0}, for each $P\in \mathcal{S}$, there exists a word
$W$ such that $A/P\cong A_{W^{\omega}}$.  If $\mathcal{S}$ is infinite, we have
infinitely many words $W_1,W_2,\ldots $ such that $W_i^j$ has nonzero image
in $A$ for every $j\ge 1$ and $A_{W_i^{\omega}}\not = A_{W_j^{\omega}}$ if $i\not =j$.  Pick integers $d_{i,j}$ such that $W_i^{d_{i,j}}$
is not a subword of $W_j^{\omega}$.
Pick an integer $m>2C$ and let $$D\ =\ \max_{1\le i,j\le m} d_{i,j}.$$
Let 
$$n \ >  \  D\max_{i\leq m}\, {\rm length}(W_{i})+D$$ be a positive integer.
By Lemma \ref{lem: 1}, there exist right infinite words $U_1,U_2,\ldots $ such that
each subword of $W_i^nU_i$ has nonzero image in $A$ and $U_i$ does not have
$W_i$ as an initial subword for $i\le m$.
Consider
the words $$Y_i = W_i^{\lfloor n/{\rm length}(W_i)\rfloor +D+1}U_i$$
for $1\le i\le m$.
Given a word $W$ with at least $b$ letters,
we let $W(a,b)$ denote the subword of $W$ formed by taking
the word whose initial position occurs at the $a$'th spot and whose final
position occurs at the $b$'th spot.
Consider the set of words 
$$\{ Y_i(j,j+n-1)~|~1\le i\le m, D{\rm length}(W_i)\le j\le n \}.$$
These words all have nonzero image in $A$.
We claim that the words in this set are distinct.
Notice that if $Y_i(k,k+n-1)=Y_j(\ell,\ell+n-1)$ then since the first $D\cdot {\rm length}(W_i)$
letters of $Y_i(k,k+n-1)$ is a subword of $W_i^{\omega}$ and the first $D\cdot {\rm length}(W_j)$
letters of $Y_i(k,k+n-1)$ is a subword of $W_j^{\omega}$, we see that
$i=j$ by definition of $D$.
But this can only occur if $k=\ell$ by definition of the words
$U_1,\ldots ,U_m$.  Thus there are at least
$$\left(n-D\max_i \, {\rm length}(W_i)\right)m$$ words of length $n$ for all $n$ sufficiently large.  Consequently, there are 
at least
$$\sum_{i=0}^n im + {\rm O}(n) \ = \ m{n\choose 2} + {\rm O}(n)$$ words of length at most $n$. 
Dividing this by $n^2$ and taking the limit as $n\rightarrow\infty$ and then using the quadratic growth hypothesis, we see 
$$m/2 \ \le \ C$$ a contradiction.  \qed
\begin{prop} Let $k$ be a field and let $A$ be a prime monomial $k$-algebra of quadratic growth.  Then
every prime homomorphic image of $A$ of GK dimension $1$ is also a monomial
algebra.\label{prop: 2}
\end{prop}
Let $N$ be the intersection of all prime monomial ideals $P$ in $A$
such that $A/P$ has GK dimension $1$.  (If $A$ has no such prime ideals, we take $N$ to be $A$.)
Note that $N\not =(0)$ since we
are taking a finite intersection of nonzero ideals and $A$ is prime.
Observe that any word in $N$ is nilpotent.  To see this, suppose that there is some word $W\in N$ such that each
subword of $W^{\omega}$ has nonzero image in $A$.  Then $A_{W^{\omega}}$ is a homomorphic image of $A$ of GK dimension $1$.  By Lemma \ref{lem: 0}, there exists some
prime monomial ideal $Q$ such that $A/Q \cong A_{W^{\omega}}$.  Notice that $W$ has nonzero image in $Q$.  But this contradicts the
fact that $W\in N\subseteq Q$.  Let $P$ be a prime ideal such that $A/P$ has
GK dimension $1$.  If $P$ is not a monomial ideal, then $N\not\subseteq P$.
Let $\overline{N}$ denote the image of $N$ in $A/P$.  We now show that
$\overline{N}$ is an algebraic ideal; that is, every element of $\overline{N}$ is algebraic over $k$.  Suppose that $\overline{N}$ is not
algebraic.  Then there is a non-algebraic element $x\in \overline{N}$; this element
is the image of a linear combination of words
$W_1,\ldots ,W_d$, all of whose images are in $\overline{N}$.  By the 
remarks above, every element in the semigroup generated by the images of
$W_1,\ldots ,W_d$ is nilpotent.  It follows from Shirshov's theorem \cite{AmS}, that
the subalgebra of $A/P$ generated by the images of $W_1,\ldots ,W_d$ is finite dimensional as a $k$-vector space.  It follows that $x$ is algebraic since it lies in this subalgebra.  Thus $\overline{N}$ is 
an algebraic ideal.  But $A/P$ is a prime Goldie ring of GK dimension $1$ and $\overline{N}$ is nonzero in $A/P$ and hence $A/(P+N)$ has GK dimension $0$.
In particular, it is finite dimensional as a $k$-vector space.  It follows that $A/P$ must be algebraic as a ring.  But a finitely generated algebra of GK dimension one is PI by the Small-Warfield theorem and hence cannot be algebraic.  Thus we obtain a contradiction.  
We conclude that $P$ is
a monomial ideal. \qed 
\vskip 2mm
To complete the proof of Theorem \ref{thm: main}, we need a result of Small about graded algebras.
\begin{lem} Let $A$ be a finitely generated prime graded algebra of quadratic growth and let $Q$ be a nonzero prime ideal of $A$ such that $A/Q$ is finite dimensional.  Then either $Q$ is the the maximal homogeneous ideal of $A$ or $Q$ contains a prime $P$ such that $A/P$ has GK dimension $1$.\label{lem: small}
\end{lem}
{\bf Proof.} Let $\mathcal{S}$ denote the collection of homogeneous ideals contained in $Q$.  Since $(0)\in \mathcal{S}$, we see that $\mathcal{S}$ is non-empty.  A standard argument using Zorn's lemma shows that $\mathcal{S}$ has a maximal element $P$.  We claim that $P$ is prime.  To see this, suppose there are nonzero $a$ and $b$ in $A$ such that $aAb\subseteq P$.  Since $P$ is homogeneous, it is no loss of generality to assume that $a$ and $b$ are homogeneous elements of $A$ (we can replace $a$ by the nonzero homogeneous part of $a$ of highest weight and do the same for $b$).   But this says that
$(AaA+P)(AbA+P)\subseteq Q$ and since $Q$ is prime, we conclude that either $AaA+P$ or $AbA+P$ is contained in $Q$, contradicting the maximality of $P$.  

We next claim that $A/P$ is PI, using an argument of Small.  Suppose that $A/P$ is not PI.
Note $A/Q$ is PI and hence satisfies a multilinear homogeneous identity; futher, by assumption $A/P$ does not satisfy this
identity.  Let $f(x_1,\ldots, x_d)$ denote this identity. 
Since $f$ is multilinear and homogeneous, there exist homogeneous elements $a_1,\ldots ,a_d$
in $A$ such that $b:=f(a_1,\ldots ,a_d)\not \in P$.  By construction $b$ is a nonzero homogeneous 
element of $Q$ that does not lie in $P$ and so $P+AbA$ is a homogeneous ideal contained in $Q$ that properly contains $P$, contradicting the maximality of $P$.  Thus $A/P$ is PI.

It remains to show that $B:=A/P$ has GK dimension $1$ if $Q$ is not the maximal homogeneous ideal of $A$.
If $Q$ is not the maximal homogeneous ideal, then $B$ is infinite dimensional and thus has at GK dimension at least $1$.  Therefore it is sufficient to show that $B$ cannot have GK dimension $2$.  We now consider $Q$ as a prime ideal of $B$ of finite codimension.  Note that $B$ is a graded prime PI algebra and hence is Goldie \cite[Corollary 13.6.6]{MR}.  It follows that if $S$ is the set of regular homogeneous elements of $B$, then we can invert the elements of $S$ to obtain an algebra $S^{-1}B =R[x,x^{-1};\sigma]$, where $R$ is a simple PI algebra \cite{GS}.  By assumption, $Q$ does not contain any nonzero homogeneous elements and hence $S^{-1}Q$ is a proper ideal of $S^{-1}B$.  

We note that $S^{-1}B/S^{-1}Q$ is finite dimensional since $B/Q$ is.  To see this, let $m$ denote the dimension of $B/Q$ and suppose $x_1,\ldots ,x_{m+1}$ are elements of $S^{-1}B$.  Then there exists some $s\in S$ such that
$sx_1,\ldots ,sx_{m+1}\in B$.  Since $B/Q$ is $m$ dimensional, some linear combination of
$sx_1,\ldots ,sx_{m+1}$ lies in $Q$.  But this says that some linear combination of $x_1,\ldots ,x_{m+1}$ lies in $S^{-1}Q$.  Hence any $m+1$ elements of $S^{-1}B$ are linearly dependent mod $S^{-1}Q$ and so $S^{-1}B/S^{-1}Q$ is finite dimensional.

Since $R$ is simple, it must embed in any homomorphic image of $S^{-1}B=R[x,x^{-1};\sigma]$ and in particular $R$ embeds in $S^{-1}B/S^{-1}Q$.  Thus $R$ is finite dimensional by the remarks above.  It follows that $S^{-1}B=R[x,x^{-1};\sigma]$ has GK dimension exactly $1$ since $R$ has GK dimension $0$.  Thus $A/P$ has GK dimension $1$.  This completes the proof.  \qed

\noindent {\bf Proof of Theorem \ref{thm: main}:}  This follows immediately from Propositions \ref{prop: 1} and \ref{prop: 2}. \qed
\section{Graded algebras of quadratic growth}
In this section, we prove Theorems \ref{thm: asmain1} and \ref{thm: asmain2} and obtain our results about primitivity for prime monomial algebras of quadratic growth.  Small \cite[Question 3.2]{Bell} asks whether a prime affine Noetherian algebra of GK dimension $2$ is either primitive or PI.  In fact, the additional hypothesis that the algebra be semiprimitive is probably necessary over countable fields to obtain this result.  We show that if $A$ is a semiprimitive prime affine monomial algebra of quadratic growth then $A$ is primitive and, moreover, it has bounded matrix images.

\begin{lem} \label{lem: tran} Let $K$ be a field, let
$A=\bigoplus_{i=0}^{\infty}A_{i}$ be a finitely generated graded
prime $K$-algebra, and let $Z$ denote the extended centre of $A$.
Suppose that $I$ is an ideal in $A$ that does not contain a nonzero homogeneous element and $z\in Z$, $x,y\in A$ are such that:
\begin{enumerate}
\item $x$ is a nonzero homogeneous element;
\item $y$ is a sum of homogeneous elements of degree smaller that the degree of $x$;
\item $x+y\in I$;
\item $zx=y$.
\end{enumerate}
Then $z $ is not algebraic over $K$.
 \end{lem}
{\bf Proof.}
 Suppose that $z$ is algebraic over $K$.  We note that $z\not =-1$ since $(1+z)a=b\not =0$.  
 Since the extended centre of a prime ring is a field (cf. Beidar et al. \cite[p. 70]{bmm}), 
 we have $(1+z)p(z)=1$ for some polynomial $p(z)$ of degree $d$.  Note that for $j\le d$ we have
 $$z^j (a_1 x a_2 x \cdots x a_d x) \ = \ a_1 y a_2 y \cdots a_j y a_{j+1} x \cdots a_d x.$$
 In particular, $p(z)(a_1 x a_2 \cdots x a_d)\in A$ for all $a_1,\ldots ,a_d\in A$.  
 Note that \begin{eqnarray*}
 xa_1 xa_2 \cdots x a_d x &=& p(z)(1+z) xa_1 xa_2 \cdots x a_d x\\
 &=& (x+y)p(z)xa_1\cdots xa_d x \\
 &\in & (x+y)A \\
 &\subseteq & I. \end{eqnarray*}
 By assumption, $I$ does not contain any nonzero homogeneous elements and hence
 $$xa_1 x a_2\cdots x a_d x \ = \ 0$$ for all homogeneous elements $a_1,\ldots ,a_d\in A$.  It follows that
 $(x)^{d+1}=(0)$.  This is impossible since
 $A$ is a prime algebra. This completes the proof. \qed
\begin{cor} Let $K$ be a  field and let $A=\bigoplus_{i=1}^{\infty }A_{i}$ be a finitely generated
 prime graded non-PI $K$-algebra of quadratic
 growth. If $P$ is  a nonzero prime ideal of $A$ then either $P$ is homogeneous or $A/P$ is PI.
 \label{cor: homogeneous}
 \end{cor}
{\bf Proof.} Let $P$ be a nonzero prime ideal $P$ such that $A/P$ is not PI.  We show that $P$ must be homogeneous.
Suppose not.  Let $Q$ be maximal homogeneous ideal of $A$ contained
in $P$. Then $Q$ is prime. Since $A/Q$ is not $PI$, by
the Small-Warfield theorem and Bergman's gap theorem $A/Q$ has
quadratic growth \cite{SW, KL}.
 The ideal $P+Q$ is a nonzero prime ideal in the graded algebra $A/Q$ and $P+Q$
 does not contain a homogeneous element. Let $\sum_{i=1}^{n}a_{i}\in
P+Q$ with $a_{i}$ of degree $i$ and $a_n\not =0$.
 We may assume that the number of nonzero elements $a_{i}$
 in this sum is minimal.
  Then by minimality of the number of nonzero terms we have
  $$a_{i}aa_{j}=a_{j}aa_{i} \qquad {\rm for~}i, j\leq n,~{\rm and}~a\in A.$$
  Beidar et al. \cite[Theorem 2.3.4]{bmm} show that
for each $i<n$, we have $$z_i a_n \ = \ a_i \qquad {\rm for~some}~z_i\in Z.$$
Since $A$ has quadratic growth, the extended centre of $A$ is algebraic over $K$ \cite[Theorem 1.3]{bs}.
Lemma \ref{lem: tran} then shows $a_{i}=0$ for $i<n$. Hence the image of $P$ in $A/Q$ contains the
 nonzero homogeneous element $a_{n}$, a contradiction.  Thus $P=Q$ and so $P$ is homogeneous.
 \qed
 \vskip 2mm
\noindent {\bf Proof of Theorem \ref{thm: asmain1}.}
 We let $K(X)$ be the field of rational functions in infinitely many
 commuting indeterminates $x_{i,j}$ where $i,j\in \mathbb{Z}$, and we let
 $A(X) = A\otimes_K K(X)$.
Define
\begin{equation} c_{i}=\sum_{j=1}^{n}x_{i,j}a_{j}\in A(X)\qquad {\rm  for~} i\ge 1\end{equation}
 and
\begin{equation} d_{i}=\sum_{j=1}^{n}x_{-i,j}a_{j}\in A(X)\qquad {\rm for}~ i\ge 1.
\end{equation}
We have two cases to consider.
\vskip 2mm 
\noindent {\bf CASE I:} the set $S=\{c_1c_{2}\cdots c_i u d_jd_{j-1}\cdots d_1~|~i,j\ge 0\}$ is linearly independent over $K(X)$.
\vskip 2mm
\noindent In this case it is easy to see
$${\rm dim} \left( AuA\cap \bigoplus_{i=0}^n A_i \right) \ \ge \ {n-m\choose 2},$$ where
$m$ is the degree of $u$.  
\vskip 2mm
\noindent {\bf CASE II:} the set $S=\{c_1c_{2}\cdots c_i u d_jd_{j-1}\cdots d_1~|~i,j\ge 0\}$ is linearly dependent over $K(X)$.  
\vskip 2mm
\noindent
Since $A(X)$ is a graded $K(X)$-algebra and $S$ is linearly dependent over $K(X)$, there is in fact some natural number $m$ such that the elements $c_{1}c_{2}\ldots c_{i}u d_{j}d_{j-1}\ldots
d_{1}$ with $i+j=m$ are linearly dependent over $K(X)$. Hence there is some
$k$ with $0\le k< m$ such that
\begin{equation} c_{1}c_{2}\ldots c_{k}ud_{m-k}d_{m-k-1}\cdots d_{1}\in c_{1}\cdots
c_{k}\sum_{j=1}^{m-k}K(X)c_{k+1}\cdots c_{k+j}ud_{m-k-j}\cdots
d_{1}. 
\label{eq: inc}
\end{equation}
  (We note that when $k=0$, we take the empty product $c_1\cdots c_k$ to be $1$.)
  From item (\ref{eq: inc}), we see
\begin{equation} A_{m}u A_{n}\subseteq \sum_{i=k}^{m-1}A_{i}uA_{m+n-i}\qquad {\rm for~}n> m. \label{eq: inc2}
\end{equation}
By comparing the left and
right hand sides of equation (\ref{eq: inc}) we also see
\begin{equation} (c_{1}c_{2}\ldots c_{k})uK(X)A_{m-k}\subseteq (c_{1}c_{2}\ldots c_{k})
   \sum_{i=1}^{m-i} K(X)c_{k+1}\ldots c_{i+k}uA_{m-i-k}.\label{eq: inc3} \end{equation}
(These are both consequences of the fact that we can clear denominators and substitute any values of $K$ for
 the variable coefficients $x_{i,j}$ that occur in the $d_j$ or $c_i$ and in this way obtain a spanning set for $A_j$ 
 from $d_j\cdots d_1$ or $c_1\cdots c_j$.)
 
Multiplying both sides of equation (\ref{eq: inc3}) on the left by $c_m$, we see
$$c_mc_1\cdots c_{k-1}c_k u K[X] A_{m-k} \subseteq  (c_mc_{1}c_{2}\ldots c_{k})
   \sum_{i=1}^{m-k} K(X)c_{k+1}\ldots c_{i+k}uA_{m-i-k}.$$
 By applying a $K$-automorphism of $K(X)$ that fixes each $d_i$ and sends $c_m\mapsto c_1$ and $c_j\mapsto c_{j+1}$ for $j\not = m$ and using equation (\ref{eq: inc}), we see
  $$c_1\cdots c_{k}uK[X]A_{m-i+k} \subseteq (c_1\cdots c_kc_{k+1})\sum_{i=1}^{m-k} K(X)c_{k+2}\ldots c_{i+k+1}uA_{m-i-k},$$
   where we take $c_{m+1}=c_1$.
 It follows that
  $${\rm dim}_{K(X)} c_1\cdots c_{k}uK(X)A_{m-k+1} 
  \ \le \ \sum_{i=1}^{m-i} {\rm dim}(A_{m-i-k}).$$
 Continuing in this manner we see that for any natural number $n\ge m-k$,
$$
 {\rm dim}_{K(X)} (c_1\cdots c_k K(X) u A_n)  \ \le \ \sum_{i=1}^{m-k} {\rm dim}(A_{m-i-k}).$$
It follows that
 $${\rm dim}_{K}(A_kuA_n) \ \le \ {\rm dim}(A_k) \sum_{i=1}^{m-k} {\rm dim}(A_{m-i-k})\qquad {\rm for~}n\ge m-k.$$
From equation (\ref{eq: inc2}), however, we have by an easy induction argument that
$A_{p}uA_{n}\subseteq \sum_{i=k}^{m-1}A_{i}uA_{p+n-i}$ whenever $p \ge m$ and $n>m$.

 Therefore if $d$ is the degree of $u$, then
 \begin{eqnarray*}&~&
{\rm  dim}_{K}AuAA_{m+1}\cap A_{n+p+d} \\ & \le &
 \sum_{i=k}^{m-1} {\rm dim}_KA_{i}uA_{n+p-i} \\
 &\le & \sum_{i=0}^{m-1-k} {\rm dim}_K A_i (A_k u A_{n+p-i-k}) \\
 &\le & \sum_{i=1}^m {\rm dim}_K(A_i) {\rm dim}_K(A_k) 
({\rm dim}_{K}(A_{1}+A_{2}+\cdots +A_{m-k})).\end{eqnarray*}
It follows that there is some constant $C_1>0$ such that
$${\rm dim}_K(AuAA_{m+1}\cap A_n)\le C_1,$$ for all $n\ge 0$.  A similar argument, replacing $A$ by $A^{\rm op}$ shows that there is a constant $C_2>0$ such that
$${\rm dim}_K(A_{m+1}AuA\cap A_n)\le C_2,$$ for all $n\ge 0$.  Since 
$$AuA \cap A_n \subseteq \left(A_{m+1}AuA+AuAA_{m+1}\right)\cap A_n$$ for all $n>d+2m+2$, we 
see that in this case
$${\rm dim}_K\left(AuA\cap \bigoplus_{i=0}^n A_i\right) \ \le \ Cn,$$ for some positive constant $C$. 
Finally, we note that if $A$ is prime then ${\rm dim}_K\left(AuA\cap \bigoplus_{i=0}^n A_i\right)$ grows at least quadratically with $n$ \cite[Lemma 2.3]{bs}, and so this case cannot occur if $A$ is prime.  This completes the proof.
\qed
 \vskip 2mm
 We are now ready to show that the intersection of the nonzero prime ideals $P$ such that $A/P$ has GK dimension $2$ is either empty of nonzero in the case that $A$ is a prime graded algebra of quadratic growth.
 \vskip 2mm
 \noindent
 {\bf Proof of Theorem \ref{thm: asmain2}.} 
 Let $$\alpha \ = \ \inf_{I} ~\liminf_{n\rightarrow\infty} n^{-2}  {\rm dim}_{K} \left(I\cap \bigoplus_{i=0}^n A_i\right),$$ where the infimum is taken over all nonzero homogeneous ideals $I$ in $A$.  Since $A$ has quadratic growth, we see that $\alpha\ge 1/2$.  It follows that we can find a homogeneous element $x\in A$ such that 
 $$\liminf_{n\rightarrow\infty} n^{-2}  {\rm dim}_{K} \left((x)\cap \bigoplus_{i=0}^n A_i\right) \ \le \ \alpha+1/4.$$
 Let $P$ be a nonzero ideal in $A$ such that $A/P$ has GK dimension $2$.  If $A/P$ is PI, then $P$ contains a homogeneous element since $A$ is not PI, using the same argument as in Lemma \ref{lem: small}; if $A/P$ is not PI, then $P$ is homogeneous by Corollary \ref{cor: homogeneous}.  Either way, we see that if $Q$ is the maximal homogeneous ideal contained in $P$, then $Q$ is a nonzero prime ideal in $A$.  
We claim that $x\in Q$.  To see this, let $B=\bigoplus_{i=0}^{\infty}B_i =A/Q,$ and let $\overline{I}=I+Q$, be the image of $I$ in $B$ and suppose that $\overline{I}$ is nonzero.
Then
$${\rm dim}_K \left(I\cap \bigoplus_{i=0}^n A_i\right)\  \ge \  {\rm dim}_K \left((I\cap Q)\cap \bigoplus_{i=0}^n A_i\right)+{\rm dim}_K\left(\overline{I}\cap \bigoplus_{i=0}^n B_i\right).$$
By Theorem \ref{thm: asmain1}, we see that
$$\liminf_{n\rightarrow\infty} \frac{1}{n^2} {\rm dim}_K \left(\overline{I}\cap \bigoplus_{i=0}^n B_i\right) \ge 1/2.$$
Moreover, $I\cap Q\not =(0)$ since $A$ is prime and hence
\begin{eqnarray*} &~&
\alpha+1/4 \\ &\ge & \liminf_{n\rightarrow\infty} \frac{1}{n^2} {\rm dim}_K\left(I\cap \bigoplus_{i=0}^n A_i\right)\\
&\ge &  \liminf_{n\rightarrow\infty} \frac{1}{n^2} {\rm dim}_K \left(I\cap Q\cap \bigoplus_{i=0}^n A_i\right)+
 \liminf_{n\rightarrow\infty} \frac{1}{n^2} {\rm dim}_K \left(\overline{I}\cap \bigoplus_{i=0}^n B_i\right)\\
 &\ge & \alpha+1/2,
 \end{eqnarray*}
 a contradiction.  It follows that $x\in Q$ and the result follows.  
 \qed
 \vskip 2mm
 From this result, we immediately obtain our results about primitivity. \vskip 2mm
\noindent  {\bf Proof of Corollary \ref{thm: asmain3}:} To see that $A$ has bounded matrix images, let $Q$ be a nonzero prime ideal of $A$ such that ${\rm GKdim}(A/Q)$ has GK dimension $0$.  By Lemma \ref{lem: small}, if $Q$ is not the maximal homogeneous ideal, then it contains a prime ideal $P$ with ${\rm GKdim}(A/P)=1$.  By Theorem \ref{thm: main} there are only finitely many such prime ideals $P$ and hence if $I$ is the intersection of all such ideals, we see that $A/Q$ is a homomorphic image of $A/I$.  Since $A/I$ is PI, we see that $A/Q$ necessarily satisfies the same polynomial identity of $A/I$ and so the matrix images are necessarily bounded.  Let $I$ denote the intersection of the maximal homogeneous ideal and all prime ideals $P$ such that $A/P$ has GK dimension at least one.  By Theorems \ref{thm: main} and \ref{thm: asmain2}, $I$ is nonzero since $A$ is prime.    Moreover, by Lemma \ref{lem: small} every nonzero primitive ideal must contain $I$.  Thus either $(0)$ is a primitive ideal and so $A$ is primitive, or the Jacobson radical of $A$ contains $I$ and so $A$ has nonzero Jacobson radical.  We note that the Jacobson radical of a monomial algebra is locally nilpotent \cite{BF}.   This completes the proof.  \qed
\vskip 2mm
We note that there do exist examples of prime monomial algebras of quadratic growth with nonzero Jacobson radical \cite{SV}.  If an algebra simply has GK dimension 2 and is not of quadratic growth then the collection of primes $P$ such that $A/P$ has GK dimension $2$ can be very strange; for example, they need not satisfy the ascending chain condition \cite{jb}.
We conclude this section by posing the following question.
\begin{quest} (Trichotomy question) Let $A$ be a prime finitely generated monomial algebra.  Is it true that $A$ is either primitive, PI, or has nonzero Jacobson radical?
\end{quest}

\section{Construction}
We now show that if we replace quadratic growth with GK dimension $2$ in Theorem \ref{thm: main}, then the conclusion of the statement of the theorem may not hold; that is, we can find a prime monomial algebra of GK dimension $2$ with unbounded matrix images.  We begin with some definitions.
\begin{defn} \em{We say that a right infinite word on a two letter alphabet $\mathcal{X}$ is \emph{Sturmian} if it has exactly
$n+1$ subwords of length $n$.}
\end{defn}
We note that if a right infinite word has fewer than $n+1$ subwords for some $n$, then it is necessarily ultimately periodic.  Thus, Sturmian words have the smallest possible subword complexity without being eventually periodic.  
We refer the reader to Allouche and Shallit \cite{AS} for examples of such words.
\begin{defn} \em{We say that a right
 infinite word on a finite alphabet $\mathcal{X}$ 
is \emph{prime} if every subword occurs infinitely often.}
\end{defn} 
We note that any right infinite Sturmian word $W$ is necessarily prime, since if $W$ has a subword $W_1$ of length $d$ that does not occur infinitely often, then we can write $W=W'V$ where $W'$ is a word and $V$ is a right infinite word with no occurrences of $W_1$.  Then the number of subwords of $V$ of length $d$ is at most $d$, since $W_1$ is not a subword of $V$.  Consequently $V$ is ultimately periodic, and thus $W$ is also ultimately periodic, a contradiction.  We note that a doubly infinite Sturmian word need not have the property that every subword occurs infinitely often; for example, consider the word
$$\cdots xxxxyyyyy\cdots $$
and the subword $xy$.
The reason we use the word \emph{prime} to describe such words comes from the
following result.
\begin{prop} Let $W$ be a right infinite word on a finite alphabet $\mathcal{X} 
$.  Let $A_W$ denote the algebra obtained by taking the free algebra on 
$\mathcal{X}$ and taking any word which is not a subword of $W$ to be a 
relation.  Then $A_W$ is prime if and only if $W$ is prime.
\end{prop}
{\bf Proof.} Let $W_1$, $W_2$ be words with nonzero image in $A_W$.  Then $W_1$ is a
subword of $W$.  Hence $W=aW_1b$, where $a$ is a finite word and $b$ is
a right infinite word.  Since $W_2$ occurs infinitely often as a subword of $W$,
$W_2$ is a subword of $b$.  Hence $b=b'W_2b''$ for some word $b'$ and some
right infinite word $b''$.  It follows that $W_1b'W_2$ is a subword of $W$ and
so $W_1b'W_2$ has nonzero image in $A$.   It follows that $A_W$ is prime.  
Conversely if $A_W$ is not prime, then there exist nonzero 
words $W_1$ and $W_2$ such
that $W_1A_W W_2 = (0)$.  Consequently, $W_2$ never occurs after the first
occurrence of $W_1$ and so $W_2$ occurs only
finitely many times as a subword of 
$W$. \qed
\vskip 1mm
Let $W$ be a right infinite prime Sturmian word.
We pick subwords of $W$ as follows.  Let $W_1$ be the first letter of $W$.
Let $W_2$ denote an initial subword of $W$ of length at least $2$
which ends with $W_1$ (such a $W_2$ exists since $W$ is prime).  In general,
having chosen $W_1,\ldots ,W_{n-1}$, we choose $W_n$ to be an initial
subword of $W$ of length
at least $2{\rm length}(W_{n-1})$ 
that ends with $W_{n-1}$.  We now define a sequence of words
$U_1,U_2,\ldots $ as follows.
We take $U_1 = W_1$ and in general we define
$$a_{i,j} \ = \ \lceil {\rm length}(W_i)/{\rm length}(W_j)\rceil,$$

\begin{equation}
V_n \ = \ W_nW_{n-1}^{a_{n,n-1}} \cdots W_2^{a_{n,2}}W_1^{a_{n,1}}W_2^{a_{n,2}}
\cdots W_{n-1}^{a_{n,n-1}}W_n, \label{eq: Vn}
\end{equation}

and
\begin{equation}
\label{eq: Un}
U_n = \big(
U_{n-1}V_n \big)^2.
\end{equation}
By induction, we see that
\begin{eqnarray}
&~&{\rm length}(U_d) \nonumber \\
& \le & 4d{\rm length}(W_d) + 8(d-1){\rm length}(W_{d-1})
+\cdots + 2^{d+1}{\rm length}(W_1)\nonumber \\
& \le & 4d^2{\rm length}(W_d)\label{eq: 4d}.\end{eqnarray}
Notice that $U_n$ is an initial subword of $U_{n+1}$ for every $n$ and hence we can define the right infinite word \begin{equation}
\label{eq: U}
U \ : = \  \lim_{n\rightarrow \infty } U_n.
\end{equation}
\begin{lem} Let $W$ be a right infinite prime Sturmian word and let
$V_1,V_2,\ldots $ be the words defined in equation (\ref{eq: Vn}).
If $d$ and $n$ are positive integers satisfying ${\rm length}(W_d)\le n<{\rm length}(W_{d+1})$, then
there are at most $12d^2(n+1)$ words of length $n$ that occur as a subword of
$V_j$ for some $j>d$.\label{lem: est}
\end{lem}
{\bf Proof.}  We let $g(n)$ denote the number of words of length $n$ that are a subword of $
V_j$ for some $j>d$.  Consider $V_j$ for $j>d$.  If we write it as a word in $W_1,W_2,\ldots $, then
$W_i$ is adjacent to $W_j$ only if $|i-j|\le 1$.  Moreover, by construction,
it can be broken into blocks of the form $W_i^k$ where ${\rm length}(W_i^k)>n$.
Thus any subword of $V_j$ is a subword of $W_i^kW_j^{\ell}$ with $|i-j|\le 1$.
If $i,j>d$, then $W_i$ and $W_j$ have length at least $n$ and both have
$W_{d+1}$ as an initial and terminal subword.
Since $W_{d+1}$
has length at least $n$, we see that if $i,j>n$,
a subword of $W_i^kW_j^{\ell}$ of length
$n$ is either a subword of $W_i$, a subword of $W_j$, or a subword of $W_{d+1}^2$.
Since we have already accounted for such words, we may consider subwords of
$W_i^kW_j^{\ell}$ with $|i-j|\le 1$ and $i$ or $j$ less than or equal to $d$.
There are at most ${\rm length}(W_m)$ subwords of length $n$ in $W_m^p$ and
at most $n+1$ subwords formed by taking a terminal subword of $W_i^k$ and an
initial subword of $W_j^k$.  Hence there are at most
$${\rm length}(W_i)+{\rm length}(W_j)+ n+1$$ subwords of $W_i^kW_j^{\ell}$
of length $n$.
Thus by considering subwords of $V_j$ for $j>d$ we get at most
$$2d(n+1+2{\rm length}(W_n)) \ \le \ 6d(n+1)$$ unaccounted for words of length
$n$.
Thus
$$g(n) \ \le \ 6d(n+1)+2(n+1)+{\rm length}(U_d).$$
By equation (\ref{eq: 4d}), we see that
$$
{\rm length}(U_d)  \le  4d^2{\rm length}(W_d)\le 4d^2 n.$$
Thus $$g(n) \ \le \ (n+1)(2+6d+4d^2)  \ \le \ 12(n+1)d^2.$$
The result follows.
 \qed
\begin{lem} Let $n\ge 2$ and let $U$ be the word defined in equation
(\ref{eq: U}).
Then the number of subwords of $U$ of length $n$ is at most
$100(n+1)(\log_2 n)^2$.
\end{lem}
{\bf Proof.} 
Let $d$ be the largest integer
such that ${\rm length}(W_d)\le n$.
Since the length of $W_j$ is at least twice the length of 
$W_{j-1}$, we see that $$1\le d\le \lfloor \log_2 n\rfloor+1.$$
Observe that
$$U \ = \ (U_dV_{d+1}U_dV_{d+1})V_{d+2}(U_dV_{d+1}U_dV_{d+1})V_{d+2}\cdots ,$$
where $U_i$ and $V_i$ are as defined in equations (\ref{eq: Un}) and (\ref{eq: Vn}).

Consequently, any subword of $U$ of length $n$ is either:
\begin{itemize}
\item{
a subword of $V_jU_dV_k$ for some $j,k>d$; or}
\item{ 
a subword consisting of a terminal subword
of $V_j$ and an initial subword of $V_{j+1}$ for some $j\ge d$.}
\end{itemize}
Since $W_{d+1}$ is an initial and terminal subword of $V_j$ for $j>d$ and has 
length at least $n$, we see that any subword of 
$U$ of length $n$ is either a subword
of $W_{d+1}U_dW_{d+1}$, a subword of $W_{d+1}^2$ 
or a subword of $V_j$ for some $j>d$.  

There are exactly
$n+1$ subwords of $W$ of length $n$.  Any subword of $W_{d+1}^2$ which is
not a subword of $W$ must consist of a terminal subword of $W_{d+1}$ and
an initial subword of $W_{d+1}$.  Since there are at most $n+1$ such words,
there are at most $2(n+1)$ words of length $n$
 which are subwords of either $W$ or $W_{d+1}^2$.  Any subword of
$W_{d+1}U_dW_{d+1}$ which is not a subword of $W$ must contain a part of $U_d$.
Hence there are at most ${\rm length}(U_d)$ such words of length $n$.
Let $f(n)$ denote the number of subwords of $U$ of length $n$.
Then we have
$$f(n)  \le  2(n+1) + {\rm length}(U_d) + \# \,{\rm of~ subwords~ of~ length ~}
n
~{\rm of~ some~}V_j, j>d.$$
Using Lemma \ref{lem: est} and equation (\ref{eq: 4d}), we see
$$f(n)\le 2(n+1)+16d^2n.$$
Since $d\le \log_2 n+1$, we see that
$$f(n) \le 16(n+1)(\log_2 n + 1)^2 + 2(n+1) \le 100(n+1)(\log_2 n)^2$$
for $n\ge 1$. \qed
\vskip 2mm
Thus we see that the number of subwords of $U$ of length at most $n$ does not grow too fast.  Our ultimate goal is to show that the algebra $A_U$ is a prime algebra of GK dimension 2 with unbounded matrix images.  We now show that $A_U$ is prime.
\begin{lem} Let $U$ be the word defined in equation
(\ref{eq: U}).
Then $U$ is prime.\label{lem: prime1}
\end{lem}
{\bf Proof.} To see that $U$ is prime, 
let $a$ be a subword of $U$.  Suppose that $a$ only appears finitely
many times as a subword of $U$. 
Then there is some $m$ such that all occurrences
of $a$ occur inside the finite word $U_m$.  But $U_{m+1}=U_mbU_mb$ for some
word $b$ and hence $a$ appears at least once more in $U_{m+1}$ than it appears
in $U_m$.  This is a contradiction.
\qed
\vskip 2mm
We use this last result along with the estimates to describe $A_U$.
\begin{prop}\label{prop: U1} Let $U$ be the word defined in equation
(\ref{eq: U}).  Then
$A_U$ is a finitely generated prime monomial algebra of GK dimension $2$.
\end{prop}
{\bf Proof.} Let $V$ be the subspace of $A$ spanned by $1$ and the elements
of the finite alphabet $\mathcal{X}$.
Then \begin{eqnarray*} {\rm dim}(V^n) & = & 1+f(1)+\cdots +f(n)
\\
&\le & 1 + f(1)+\sum_{j=2}^n 100(j+1)(\log_2 j)^2
 \\
&\le & 1+f(1)+100(n+1)^2(\log_2 n)^2. \end{eqnarray*}
Thus for every $\varepsilon>0$ we have 
$${\rm dim}(V^n) < n^{2+\varepsilon}\qquad {\rm for}~n~{\rm sufficiently~large}.$$
It follows that $A_U$ has GK dimension at most $2$.  On the other hand,
the images of the distinct subwords of the Sturmian word $W$ are linearly independent.
Hence ${\rm dim}(V^n)\ge {n+1\choose 2}$ and so $A_U$ has GK dimension
at least $2$.  The fact that $A_U$ is prime follows from Lemma \ref{lem: prime1}.  \qed
\vskip 2mm
We now show that the algebra $A_U$ has unbounded matrix images.  To do this we need a simple estimate for the PI degree of a class of monomial algebras.
\begin{lem} Let $T$ be a right infinite periodic word with minimal period $d$.
Then $A_T$ is a prime ring of GK dimension $1$ satisfying $S_{2d}$ and satisfying
no identity of smaller degree.
\end{lem}
{\bf Proof.}
 Since $T$ is periodic, it is prime and hence $A_T$
is a prime ring.
Since $T$ is periodic and infinite, $A_T$ 
has GK dimension $1$.  If follows from the Small-Warfield theorem \cite{SW} that $A_T$
is PI.
Let $Y_1,\ldots ,Y_d$ denote the $d$ distinct words of $T$ of length $d$.  Then it is
easy to check that the image of $$Z \ := \ Y_1+\cdots +Y_d$$ is central in $A_T$.
Observe that the image of $Y_iY_j$ is $0$ if $i\not =j$ and $Y_i^2 = Y_iZ = ZY_i$ in $A_T$.  Thus 
$Y_1Z^{-1},\ldots ,Y_nZ^{-1}$ is a set of orthogonal idempotents in the 
quotient ring $Q(A_T)$ of $A_T$.  Since $Q(A_T)$ is a simple Artinian ring,
 $A$ cannot satisfy an identity of
degree less that $2d$. Let $K$ denote the centre of the quotient ring of $A_T$.
We claim that the quotient ring of $A_T$ embeds in $M_d(\overline{K})$.  
For $1\le i,j\le
d$ let 
$T_{i,j}$ denote the subword of $T$ starting at the $i^{\rm th}$ position of
$T$ and terminating at the $(2d+j)^{\rm th}$ position of $T$.  Let $t\in
\overline{K}$ be a primitive $d^{th}$ root of $Z$ and
let $e_{i,j} = T_{i,j}t^{i-j-2d}$.  Then $\{e_{i,j}\}$ is a set of
matrix units.  Moreover, any element of $A_T$ is in the $\overline{K}$-span
of these matrix units.  Hence we get the desired embedding.  It follows that
$A_T$ satisfies $S_{2d}$ and no smaller identity.  \qed
\begin{prop} Let $U$ be the word defined in equation
(\ref{eq: U}).
Then the algebra $A_U$ has unbounded matrix images.\label{prop: U2}
\end{prop}
{\bf Proof.} Observe that for each $j$ and $m$, $W_j^m$ is a subword of $U$.
Let $T_j$ denote the periodic right infinite word $W_jW_jW_j\cdots $.  Observe that we have a homomorphism
$$A_U\rightarrow A_{T_j}$$ since every subword of $T_j$ is a subword of $U$.
Let $d_j$ denote the minimal period of $T_j$.  Then $A_{T_j}$ has PI degree
$2d_j$.  To show that $A_U$ has unbounded matrix images, we must show that
$\lim_{j\rightarrow \infty} d_j = \infty$.  To see this, suppose that this is
not the case.  Then there exists some $m$ such that $d_j=m$ for infinitely many $j$.  In particular, $W_j$ is periodic with period $m$ for infinitely many $j$.
Since $W_j$ is an initial subword of $W$ for every $j$ and 
${\rm length}(W_j)\rightarrow\infty$, we see that $W$ must also be periodic and have period at most $m$, contradicting the fact that $W$ has at least $n+1$ subwords of length $n$ for every $n$. 
\qed  
\vskip 2mm
\noindent {\bf Proof of Theorem \ref{thm: main2}:} This follows easily from Proposition \ref{prop: U1} and \ref{prop: U2}. \qed
\section*{Acknowledgments}
We thank Tom Lenagan for many helpful comments and suggestions.

\end{document}